\documentclass[12 pt ]{article}

\usepackage{amsmath,amsthm,amsfonts}

\usepackage[T2A]{fontenc}
\usepackage[cp1251]{inputenc}
\usepackage[english]{babel}
\usepackage{amsmath}
\usepackage{amssymb}
\usepackage{longtable}

\allowdisplaybreaks[1]

\newtheorem{thm}{Theorem}

\everymath{\displaystyle}
\theoremstyle{remark}

\theoremstyle{definition}

\begin{document}

\title{The proof of conjecture of Brutman and Passow \begin{footnote}{The research was partly supported by the Sofia University Research Foundation
 through Contract no 139/2016}\end{footnote}}
 
\author{Milko Damyanov Takev}


\maketitle

\newtheorem{dd}{Theorem}
\newtheorem{cc}{Lemma}
\begin{abstract}

We consider the conjecture  of Brutman and Pasow on a totality divided differences
 and prove the conjecture for continuous functions.

\end{abstract}

\textit{ \textbf {Key words:}} divided differences problem, finite-differences.

\bigskip
Let $C[0,1]$ be the set of continuous
 functions, defined on the segment  $[0,1]$and equipped with the uniform norm $ \vert\vert .\vert\vert $.
 We denote by $\epsilon $ any infinitely small positive number.
Divided difference of $f$ of order $n$ at the distinct knots $x_{i}, i =0,1\ldots , n$ is

$$ f[x_{0},\ldots ,x_{n}] :=\sum_{k=0}^{n} {f(x_{k})\over  \Omega^{'}(x_k)},$$
where $\Omega(x):= (x-x_{0})\ldots (x-x_n)$.
As usual we define the finite differences:
\begin{eqnarray*}
\Delta_{h}^{n}f(x) &=&
\sum_{j=0}^{n}
(-1)^{n+j} 
f(x+jh ).\\
\end{eqnarray*}

From binomial theorem  it is easy to see that:
\begin{eqnarray}\label{tem}
f[0,1/n,2/n,...,1] = {\Delta_{1\over n}^{n}f(0)\over n!({1\over n}    )^{n}}
\end{eqnarray}

 In \cite{kn: br} is 
formulated the following:

\it{
Conjecture (Brutman and Pasow).
Let $f \in C[0,1]$. Suppose that the divided differences
$f[0,1/n,2/n,...,1] = 0, n=1,2,... $ and $f(0)=0,$
then $f(x) = 0$ for all $x\in [0,1]$.}
\rm

Let $\Delta_{1\over n}^{n}f(0)= 0$ and
$ Q_{n}(f;x)$ be the Interpolation polynomial of Lagrange of  degree $n$ at the knots
$ 0,1/n,2/n,...,1$. It is  known that the leading coefficient of
$ Q_{n}(f;x)$ is $f[0,1/n,2/n,...,1]$. Using this,
we may formulate the conjecture of Brutman and Pasow in the following equivalent form:

\it{
Let $f \in C[0,1]$  and $f(0)=0$ suppose that 
the Interpolation polynomial of Lagrange of  degree $n$ at the knots
$ 0,1/n,2/n,...,1$. is in fact of  degree $n-1$ for each $n \geq 1$.
Does this imply that$f(x) = 0$ for all $x\in [0,1]$.}
\rm{

Let us mention, that the conjecture of Brutman and Pasow is true for polynomials. 
Really, let $f$ be polynomial
$ f = a_nx^n + a_{n-1}x^{n-1} + \ldots + a_0$ and from integral representation of
finite-differences
$$\Delta_{h}^{n}f(x) = \int_0^h \ldots \int_0^h f^{(n)}(x+t_1 + \ldots  +t_n) dt_1 \ldots dt_n$$
we in strict succession obtain, that
$a_n = a_{n-1}= \ldots = a_0 = 0 $ and $f\equiv 0.$
  In {\cite{kn: mt}} the conjecture was proved  for entire  functions of 
 exponential type.

\thm{\label{te.S}
Let $f\in C[0,1]$ for any fixed $\rho$ and $f(0)=0$.
Suppose that the
 interpolation polynomial of Lagrange
 of  degree $n$ at the knots
$ 0,1/n,2/n,...,1$ is in fact of  degree $n-1$ for each $n \geq 1$.
Then $f(x) = 0$ for all $x\in [0,1]$}.
}
\rm
\\ Proof. 
First, we set 
for any $ x \notin [0,1], f(x)=0.$ Take any $x$ and $y$ in $[0,1]$.
Say $y = x +n!h$. Now $x<y$ for each integer $j=1,\ldots , n,$ write $n! = jj'.$. We use
the formula for infinity interval from \cite{kn: hw} and obtain

\begin{eqnarray*}
(-1)^{n}[f(x) -  Q_{n}(f;x) -f(y) -  Q_{n}(f;y)] \\=
{1\over {N}}\sum_{i=1}^{N} \{\Delta_{ih}^{n}[f(x)-Q_{n}(f;x)] - \Delta_{ih}^{n}[f(y)-Q_{n}(f;y)]\}\\
-{1\over {N}} \sum_{j=1}^{n} (-1)^{n+j}{n\choose j} \sum_{i=1}^{N}[f(x)-Q_{n}(f;x) - f(y)+Q_{n}(f;y)]
\end{eqnarray*}
Now we put $y=1,$ and have $f(y+ijh) = 0$ for every $ij,Q_{n}(f;y)= 0.$ For the first term ot right side

\begin{eqnarray}
{1\over {N}}\sum_{i=1}^{N} \{\Delta_{ih}^{n}[f(x)-Q_{n}(f;x)] - \Delta_{ih}^{n}[f(y)-Q_{n}(f;y)]\}
= (-1)^{n}f(x) + \epsilon
\end{eqnarray}

\begin{eqnarray*}
 {1\over {N}} \sum_{j=1}^{n} (-1)^{n+j}{n\choose j} \sum_{i=1}^{N}[f(x)-Q_{n}(f;x) - f(y)+Q_{n}(f;y)]\\ ={1\over {N}}
 \sum_{j=1}^{n} (-1)^{n+j}{n\choose j}\{ \sum_{i=N-j'+1}^{N}[f(x + jih)-Q_{n}(f;x+jih)] \\ - \sum_{i=1}^{j'}[f(y+ijh)-Q_{n}(f;y+ijh)]\}\\ =
{1\over {N}} \sum_{j=1}^{n} (-1)^{n+j}{n\choose j} \sum_{i=N-j'+1}^{N}f(x + jih)-\\
{1\over {N}} \sum_{i=N-j'+1}^{N}  \{\sum_{j=0}^{n} (-1)^{n+j}{n\choose j} Q_{n}(f;x+jih) + (-1)^{n}Q_{n}(f;x) \} \\ - 
{1\over {N}} \sum_{j=1}^{n} (-1)^{n+j}{n\choose j} \sum_{i=1}^{j'}f(y+ijh)-\\
{1\over {N}} \sum_{i=1}^{j'}  \{\sum_{j=0}^{n}(-1)^{n+j}{n\choose j}Q_{n}(f;y+ijh)+(-1)^{n}{{j'}\over {N}}Q_{n}(f;y)\} \\ 
={1\over {N}} \sum_{j=1}^{n} (-1)^{n+j}{n\choose j} \sum_{i=N-j'+1}^{N}f(x + jih)-  \\ 
{1\over {N}} \sum_{j=1}^{n} (-1)^{n+j}{n\choose j}\sum_{i=1}^{j'}f(y+ijh)+(-1)^{n}{{j'}\over {N}}Q_{n}(f;y)\\
= (-1)^{n}{{j'}\over {N}}Q_{n}(f;x) + \epsilon.
\end{eqnarray*}
Using (2) and above considaration we get
\begin{eqnarray}
(-1)^{n}[f(x) -  Q_{n}(f;x) ] =(-1)^{n}f(x) - (-1)^{n}{{j'}\over {N}}Q_{n}(f;x) + \epsilon.
\end{eqnarray}

From (3) we obtain that $Q_{n}(f;x) \equiv 0$. We obtain $f(x) = 0$ in interpolation knots $0, {1\over n}\ldots {n\over n},  $
 for every natural $n$.
 Finally continuous function  $f(x) = 0$ in rational points in $[0,1]$, therefore $f \equiv o.$

 {\em Section of Mathematical Analysis, Department of Mathematics and
Informatics

Sofia University, 5 J. Bourchier, 1164 Sofia, Bulgaria}

e-mail: takev@fmi.uni-sofia.bg

  \end{document}